\documentclass[12pt]{amsart}

\usepackage{amsmath}
\usepackage{amsthm}
\usepackage{amssymb}

\usepackage{ifthen}

\input diagrams

\newcommand{\trt}{3}
\newcommand{\lbl}[1]
{\ifthenelse{\trt=2}{\quad\label{#1}\ref{#1}}{\label{#1}}}

\newcommand{\bff}{\bf}

\newfont{\rmm}{cmr10 scaled 1000}
\newfont{\itt}{cmsl10 scaled 1000}
\newfont{\nazad}{cmff10 scaled 1000}

\theoremstyle{plain}

\newtheorem{theo}{Theorem}[section]
\newtheorem{lemm}[theo]{Lemma}
\newtheorem{prop}[theo]{Proposition}
\newtheorem{coro}[theo]{Corollary}

\theoremstyle{definition}

\newtheorem{defi}[theo]{Definition}
\newtheorem{rema}[theo]{Remark}

\newcommand{\armaz}
{M.Artin, B.Mazur,
\emph{On periodic points},
Annals of Math.
{\bff 102}
(1965),
    82 -- 99.
}

\newcommand{\baladi}
{V.Baladi
\emph{Periodic orbits and dynamical spectra},
Ergodic theory and dynamical spectra,
{\bff 18}
(1998),
    255 - 292.
}

\newcommand{\fried}{  D.Fried,
\emph{Homological Identities for closed orbits},
Inv. Math.,
{\bff 71},
(1983)
419--442.
}

\newcommand{\fel}{  A.Fel'shtyn,
\emph{
Dynamical zeta functions, 
Nielsen Theory and Reidemeister torsion,}
Memoirs of AMS {\bff 147} (2000).
}

\newcommand{\friedtwi}{  D.Fried,
\emph{Periodic points and twisted coefficients},
Lect. Notes in Math.,
{\bff 1007},
(1983)
261--293.
}

\newcommand{\geo}{  R.Geoghegan, A.Nicas,
\emph{Trace and torsion in the theory of flows},
Topology,
{\bff 33},
(1994)
683 -- 719
}

\newcommand{\grayson}{  D.Grayson,
\emph{Higher Algebraic K-theory 2},
Lecture Notes in Math.,
{\bff 551}, Springer,
(1976)
217 -- 240
}

\newcommand{\jiang}{  B.Jiang,
\emph{ Estimation of the number of periodic orbits}
Pacific J. Math. {\bff 172}, 151--185 (1996)
}

\newcommand{\hulitop}{   M.Hutchings, Y-J.Lee
\emph{              Circle-valued Morse theory, Reidemeister torsion
and Seiberg-Witten invariants of 3-manifolds},
Topology,
{\bff 38},
(1999),
861 -- 888}

\newcommand{\milnWT}
{J.Milnor,
\emph{ Whitehead Torsion},
Bull. Amer. Math. Soc.
{\bff 72}
(1966),
358 - 426.
}

\newcommand{\patou}
{ A.V.Pazhitnov,
\emph{
On the Novikov
complex for rational Morse forms},
 Annales de la Facult\'e de Sciences de Toulouse,
{\bff 4}      (1995), 297 -- 338.
 }

\newcommand{\pasur}
{ A.V.Pajitnov,
\emph{
Surgery on the Novikov Complex},
K-theory
{\bff 10} (1996),      323-412.
 }

\newcommand{\pamrl}
{  A.V.Pajitnov,
\emph{Rationality and exponential growth
properties of the boundary operators in the Novikov
Complex},
Mathematical Research Letters,
{\bff 3}
(1996),
  541-548.
 }

\newcommand{\pafest}
{  A.V.Pajitnov,
\emph{Simple homotopy type of the Novikov
Complex and Lefschetz $\zeta$-function of the gradient flow},
Russian Math.Surveys, 1999, {\bf 54} No. 1 }

\newcommand{\ranpaj}
{   A.V.Pajitnov, A.Ranicki,
\emph{The Whitehead group of the Novikov ring},
to appear in "K-theory" }

\newcommand{\schutz}
{  D.Sch\"utz,
\emph{Gradient flows of closed 1-forms and their closed orbits},
math.DG/0009055}

\newcommand{\sieben}
{L.Siebenmann,
\emph{A total Whitehead torsion obstruction
 to fibering over the circle},
Comment. Math. Helv. {\bff 45} (1970), 1--48.}

\newcommand{\smdyn}
{  S.~Smale,
\emph{Differential dynamical systems},
 Bull. Amer. Math. Soc. {\bff 73} (1967)
747--817.}

\newcommand{\bour}
{N.Bourbaki,
\emph{Groupes de Lie, Alg\`ebres de Lie
}, Springer, 1985
}

\newcommand{\dold}
{   A.Dold
\emph{Lectures on Algebraic Topology
},
Springer,  1972.
}

\newcommand{\kiang}
{   Kiang Tsai-han
\emph{The theory of Fixed point classes},
 Springer, 1989.
}

\newcommand{\milnhcob}
{   J.~Milnor,
\emph{Lectures on the
h-cobordism theorem},
 Princeton University
Press,
 1965.
}

\newcommand{\ranibook}
{ A.Ranicki,
\emph{High-dimensional knot theory},
Springer Mathematical Monograph, Springer
(1998)}

\newcommand{\farran}{ M.Farber, A.Ranicki,
\emph{ The Morse-Novikov theory 
of circle-valued functions and non-commutative localization},
Proc. 1998 Moscow Conference for S.P.Novikov's 60th Birthday,
Proc. Steklov Inst. 225, 381--388 (1999); e-print dg-ga/9812122}

\newcommand{\lueckzeta}
{ W.L\"uck,
\emph{The Universal Functorial Lefschetz Invariant }
Fund. Math. 161, 167--215 (1999)
}

\def\paepri{ A.V.Pajitnov,
\emph{Incidence coefficients in the Novikov Complex
for Morse forms: rationality and exponential growth properties},
e-print dg-ga/9604004 20 Apr 96.}

\newcommand{\rann}
{  A.Ranicki,
      \emph{
The algebraic construction
of the Novikov complex of a circle-valued Morse function}
e-print, math.AT/9903090, 15 mar 99
 }

\def\novidok{S.P.Novikov,
\emph{Mnogoznachnye funktsii i funktsionaly.
Analog teorii Morsa}, Doklady AN SSSR
{\bf 260}
(1981),
 31-35.
   English translation:
 S.P.Novikov. \emph{ Many-valued functions
 and functionals. An analogue of Morse theory},
 Sov.Math.Dokl.
{\bff 24}
(1981),
       222-226. }

\newcommand{\pastpet}
{ A.Pajitnov
\emph{
      Ratsional'nost' granichnyh operatorov 
v komplekse Novikova v sluchae obschego
polozheniya}
Algebra i Analiz,
{\bff 9}, no.5 (1997), p. 92--139.
\quad
English translation:
 \emph{
      The incidence coefficients in the Novikov complex
are generically rational functions}
St. Petersburg Mathematical Journal,
{\bff 9}, no.5 (1998), p. 969--1006 }

\newcommand{\postnlie}
{M.M.Postnikov,
\emph{ Gruppy i algebry Li}
Moskva, Nauka, 1982
(in Russian)}\renewcommand{\a}{\alpha}
\renewcommand{\b}{\beta}
\newcommand{\g}{\gamma}
\renewcommand{\d}{\delta}

\newcommand{\ve}{\varepsilon}
\newcommand{\z}{\zeta}

\renewcommand{\l}{\lambda}
\renewcommand{\k}{\varkappa}
\newcommand{\m}{\mu}
\newcommand{\n}{\nu}

\renewcommand{\o}{\omega}

\newcommand{\G}{\Gamma}
\newcommand{\D}{\Delta}

\renewcommand{\L}{\Lambda}

\renewcommand{\AA}{{\mathcal A}}

\newcommand{\CC}{{\mathcal C}}

\newcommand{\EE}{{\mathcal E}}

\newcommand{\GG}{{\mathcal G}}
\newcommand{\HH}{{\mathcal H}}

\newcommand{\LL}{{\mathcal L}}
\newcommand{\MM}{{\mathcal M}}

\newcommand{\PP}{{\mathcal P}}
\newcommand{\QQ}{{\mathcal Q}}
\newcommand{\RR}{{\mathcal R}}

\newcommand{\UU}{{\mathcal U}}

\newcommand{\NNN}{{\mathbf{N}}}

\newcommand{\QQQ}{{\mathbf{Q}}}
\newcommand{\RRR}{{\mathbf{R}}}

\newcommand{\ZZZ}{{\mathbf{Z}}}

\newcommand{\gC}{{\mathfrak{C}}}

\newcommand{\gG}{{\mathfrak{G}}}

\newcommand{\gL}{{\mathfrak{L}}}

\newcommand{\gW}{{\mathfrak{W}}}

\renewcommand{\leq}{\leqslant}
\renewcommand{\geq}{\geqslant}

\newcommand{\Prf}{{\it Proof.\quad}}

\newcommand{\Wkr}{W^{\circ}}

\newcommand{\Log}{\text{\rm Log }}

\newcommand{\Ker}{\text{\rm Ker }}

\newcommand{\ind}{\text{\rm ind}}

\renewcommand{\Im}{\text{\rm Im }}
\newcommand{\supp}{\text{\rm supp }}
\newcommand{\Int}{\text{\rm Int }}

\newcommand{\Fix}{\text{\rm Fix }}

\newcommand{\smo}{C^{\infty}}

\newcommand{\chart}{\Phi_p:U_p\to B^n(0,r_p)}
\newcommand{\atlas}{\{\Phi_p:U_p\to B^n(0,r_p)\}_{p\in S(f)}}

\newcommand{\fcob}{f:W\to[a,b]}

\newcommand{\indl}[1]{{\scriptstyle{\text{\rm ind}\leqslant {#1}~}}}
\newcommand{\inde}[1]{{\scriptstyle{\text{\rm ind}      =   {#1}~}}}
\newcommand{\indg}[1]{{\scriptstyle{\text{\rm ind}\geqslant {#1}~}}}

\newcommand{\pr}{\partial}

\newcommand{\id}{\text{id}}

\newcommand{\st}[1]{\overset{\rightsquigarrow}{#1}}

\newcommand{\stexp}[1]{{#1}^{\rightsquigarrow}}

\newcommand{\qt}{\hfill\triangle}
\newcommand{\qs}{\hfill\square}

\newcommand{\pa}{\vskip0.1in}

\newcommand{\wi}{\widetilde}

\newcommand{\ove}{\overline}
\newcommand{\unde}{\underline}

\renewcommand{\(}{\big(}
\renewcommand{\)}{\big)}

\newcommand{\wh}{\widehat}

\newcommand{\stv}{\stexp {(-v)}}

\newcommand{\vvbsm}{V_b^{( s-1)    }}
\newcommand{\vvasm}{V_a^{( s-1)    }}
\newcommand{\vvbs}{V_b^{[s]}    }
\newcommand{\vvas}{V_a^{[ s]    }}

\newcommand{\factor}{\vvbs / \vvbsm}
\newcommand{\factora}{\vvas / \vvasm}

\newcommand{\vvksm}{\wi V_k^{( s-1)   }}
\newcommand{\vvks}{\wi V_k^{[s]}    }

\newcommand{\vvkmsm}{\wi V_{k-1}^{( s-1)   }}
\newcommand{\vvkms}{\wi V_{k-1}^{[s]}    }

\newcommand{\fac}{\vvks / \vvksm}
\newcommand{\facm}{\vvkms / \vvkmsm}

\newcommand{\vflesh}{v\!\da}
\newcommand{\vvfl}{\wi v\!\da}

\newcommand{\vbsm}{V_b^{\{\leq s-1\}}    }
\newcommand{\vasm}{V_a^{\{\leq s-1\}}    }
\newcommand{\vbs}{V_b^{\{\leq s\}}    }
\newcommand{\vas}{V_a^{\{\leq s\}}    }

\newcommand{\bere}{\begin{rema}}
\newcommand{\bede}{\begin{defi}}

\renewcommand{\beth}{\begin{theo}}
\newcommand{\bele}{\begin{lemm}}
\newcommand{\bepr}{\begin{prop}}
\newcommand{\beeq}{\begin{equation}}
\newcommand{\bega}{\begin{gather}}
\newcommand{\been}{\begin{enumerate}}

\newcommand{\beal}{\begin{aligned}}

\newcommand{\enre}{\end{rema}}

\newcommand{\enpr}{\end{prop}}
\newcommand{\enth}{\end{theo}}
\newcommand{\enle}{\end{lemm}}
\newcommand{\enen}{\end{enumerate}}
\newcommand{\enga}{\end{gather}}
\newcommand{\eneq}{\end{equation}}
\newcommand{\enal}{\end{aligned}}

\newcommand{\subs}{\subsection}

\newcommand{\lb}{\label}

\newcommand{\Lxi}{{\wh \L}_\xi}
\newcommand{\LxiQ}{{\wh \L}_{\xi,\QQQ}}

\newcommand{\tens}[1]{\underset{#1}{\otimes}}

\newcommand{\da}{\downarrow}

\newcommand{\dow}{\pr_0 W}

\newcommand{\daw}{\pr_1 W}

\newcommand{\Tr}{{\text Tr}}

\newcommand{\limdir}{\underset {\to}{\lim}}

\newcommand{\bu}{\bullet}

\newcommand{\kom}[2]{ {#1}{#2}{ {#1}^{-1}} {{#2}^{-1}} }

\newcommand{\kommm}[2]{ {#1}'{#2}'{ ({#1}'')^{-1}} {({#2}'')^{-1}} }

\newcommand{\cmxi}{\wh C_*^\D( M, \xi)}
\newcommand{\whgxi}{\wh {{\rm Wh}} (G,\xi)}
\newcommand{\klxi}{K_1(\Lxi)}
\newcommand{\klxiQ}{K_1(\wh \L_{\xi, \QQQ})}

\newcommand{\whgxiQ}{\wh {{\rm Wh}}_\QQQ (G,\xi)}

\newcommand{\arar}{A_\rho}

\newcommand{\arz}{\arar [[t]]}

\newcommand{\arzz}{\arar ((t))}

\newcommand{\karzz}{K_1\big(\arzz\big)}

\newcommand{\karz}{K_1(\arz)}

\newcommand{\bqq}{\begin{equation*}}
\newcommand{\GL}{\text{\rm GL }}

\newcommand{\bq}{\begin{equation}}

\newcommand{\sut}{~such~that~}

\newcommand{\wrt}{~with respect to}
\newcommand{\sm}{\setminus}

\newcommand{\sbs}{\subset}
\newcommand{\ho}{homomorphism}

\newcommand{\ma}{manifold}
\newcommand{\nei}{neighborhood}
\newcommand{\dfm}{diffeomorphism}

\newcommand{
\sma}{submanifold}

\newcommand{\noconf}{~there~is~no~possibility~of~confusion}

\newcommand{\ATA}{Almost~ Transversality~ Condition}

\newcommand{\hot}{homotopy}

\newcommand{\TA}{Transversality Condition}

\newcommand{\hog}{homology}

\newcommand{\su}{subsection}

\begin{document}

\title[Gradient flows and non abelian Witt vectors]
{ Closed orbits of gradient flows \\
and logarithms of non-abelian Witt vectors}
\author{A.V.Pajitnov}
\address{
UMR 6629 CNRS, Universit\'e de Nantes\\
D\'epartement de Math\'ematiques\\
2, rue de la Houssini\`ere, 44072, Nantes Cedex France}
\email{pajitnov@math.univ-nantes.fr}
\begin{abstract}
We consider the flows generated by generic gradients
of Morse maps $f:M\to S^1$.
To each such flow
we associate an invariant counting
the closed orbits of the flow. Each closed orbit
is counted with the weight
derived from its index and homotopy class.
The resulting invariant is called the {\it eta function},
and lies in a suitable quotient of the Novikov completion
of the group ring of the fundamental group of $M$.
Its abelianization coincides with the logarithm of
the twisted Lefschetz zeta function of the flow.
For $C^0$-generic gradients
we obtain a formula expressing the eta function
in terms of the torsion of a special homotopy equivalence
between the Novikov complex of the gradient flow and the completed
simplicial chain complex of the universal cover.
\end{abstract}
\maketitle

\section*{Introduction }
\lb{s:intro}

The study of the periodic orbits of a dynamical system
via homotopy invariants
is one of traditional subjects of the algebraic topology
and dynamical system theory.
Let us consider first a discrete dynamical system
given by a continuous map
$f:X\to X$
of a topological space $X$.
The fundamental invariant of this dynamical system --
the cardinality of
the fixed point set
$\Fix f$ --
has a homotopy counterpart, namely
the algebraic number of fixed points of $f$.
The Lefschetz trace formula expresses this
number (under some mild conditions on $f$ and $X$)
in terms of the traces of the \ho s
induced by $f$ in homology.

In order to investigate the sets
$\Fix f^n$
for $n\to\infty$
Artin and Mazur
\cite{armaz}
introduced a zeta function
which encodes the information about all the periodic
points into a single power series in one variable $t$.
This zeta function and its generalizations have
proved very important in dynamical system theory.
See the works \cite{baladi} of V.Baladi
and \cite{fel} of A.Fel'shtyn
for a survey of the present state of the theory.

The Artin-Mazur zeta function
has a homotopy counterpart --
the {\it Lefschetz zeta function }
introduced by S.Smale \cite{smdyn}.
Here is the definition.
Let $\Fix f$
denote the set of
fixed points of $f$.
Assume for simplicity that for any $n$
the set
$\Fix f^n$
 is finite.
Each  $x\in \Fix f^k$ has an associated index $\nu(x)\in\ZZZ$.
Set
\bq
L_k(f)=\sum_{a\in \Fix f^k} \nu(a)
\end{equation}
and define the Lefschetz zeta function by the following formula:

\begin{equation}
\zeta_{L}(t)=\exp \bigg(\sum_{k=0}^\infty \frac {L_k(f)}k t^k\bigg)
\lbl{f:lefzet}
\end{equation}

Here is the formula for $\zeta_L$
in terms of the homology invariants of  $f$:
\begin{equation}
\zeta_{L}(t)= \prod_i \det(I-tf_{*,i})^{(-1)^{i+1}}
\lbl{f:detzet}
\end{equation}
where
$f_{*,i}$
stands for
the \ho~induced in $H_i(X)$ by $f$.

In order to encode  more information on periodic points in one
single power series, one counts  the periodic points
with weights belonging to some group ring
(usually the group ring of some regular covering of $X$).
An abelian  zeta function of this kind was first introduced by D.Fried
 \cite{friedtwi},  and the corresponding analog of the
formula (\ref{f:detzet})
was obtained.
The first non-abelian generalization
of zeta functions was introduced by R.Geoghegan and A.Nicas
\cite{geo} as one of the consequences of
their parameterized
Lefschetz-Nielsen fixed point theory.
Their
{\it Nielsen-Fuller invariant}
counting periodic points of an orientation
preserving diffeomorphism
$f:M\to M$
belongs to the Hochschild homology
of the group ring of
the mapping torus of $f$.
They obtain an analog of the formula
(\ref{f:detzet})
for this invariant.
The analogs of Lefschetz zeta functions for continuous maps
in the non-commutative setting were defined also
in the papers \cite{jiang}, \cite{lueckzeta}.

Proceeding to the dynamical systems generated by
flows on manifolds, one would expect that the theory
described above generalizes to this setting.
Namely, there should exist an analog of Lefschetz
zeta function
(\ref{f:lefzet})
counting the closed orbits of the flow, and there should
be a formula expressing this invariant
in homotopy invariant terms.
Such formulas exist in many important cases, but
only for non-singular flows.
See \cite{fried} Th.10 for abelian invariants
of homology proper Axiom A-No-Cycles semiflows
and \cite{geo} for non-abelian invariants of suspension flows.

The first results concerning zeta functions of
flows with zeros
appeared only very recently
(\cite{hulitop}, \cite{pafest}).
These papers deal with the flows generated by gradients of
 circle-valued Morse maps.
The results of these papers show that the naturally arising zeta 
functions
are no longer homotopy invariants of
the underlying manifold itself (as it is  the case for
non-singular flows).
There occurs a "correction term"
which can be computed in terms of the Novikov complex (see 
\cite{novidok})
associated with the flow.
Both papers cited above deal with abelian invariants.

In the present paper we introduce
{\it non-abelian}
invariants counting  closed orbits of
gradient flows of circle-valued Morse maps.
These invariants -- we call them {\it eta functions} --
generalize the logarithms of abelian zeta functions.
(Note that the
abelian zeta functions carry
the same amount of information as their logarithms.)
The definition of  the non-abelian eta function for gradient flows
is given in
\ref{susu:nonabeta}.
The Main Theorem of the paper
gives a formula expressing the eta function of the flow
in terms of the torsion of a chain homotopy equivalence
naturally associated with the flow.

{\it  Remarks on the terminology.}
Let $R$ be a ring with a unit.
Recall from \cite{milnWT}
the group
$\ove{K}_1 (R)=K_1 (R)/\{0, [-1]\}$
where $[-1]$ denotes  the element of order
 2 corresponding to the unit
$(-1)\in GL(R,1)\sbs GL(R,\infty)$.
A free based finitely generated chain complex $C_*$ of right 
$R$-modules
is called {\it $R$-complex}.
Recall from \cite{milnWT}, \S 3
that for an
acyclic $R$-complex $C_*$
the torsion $\tau(C_*)\in \ove{K_1}( R)$
is defined.
If $\phi:F_*\to D_*$
is a \hot~equivalence
of $R$-complexes, then its {\it torsion}
$\tau(\phi)\in \ove{K_1}(R)$
is defined to be the torsion
of the algebraic mapping cone $C_*(\phi)$.

\subsection{Counting closed orbits of a gradient flow of
a circle-valued Morse map}
\lb{su:stnov}
\pa
In order to state the main theorem  we shall need some
preliminaries on circle-valued Morse theory. They are gathered in the
three following subsections.

\subsubsection{ Novikov rings}
\label{susu:algnot}

Let
$G$ be a group, and
$\xi:G\to \ZZZ$
be a \ho.
Set $\L=\ZZZ G$
and denote by
$\wh{\wh\L}$
the abelian group of all functions $G\to\ZZZ$.
Equivalently,
$\wh{\wh\L}$
is the set of all formal linear combinations
$\l=\sum_{g\in G} n_g g$
 (not necessarily finite)
of the elements of $G$ with integral coefficients $n_g$.
For $\l\in \wh{\wh\L}$ set
$\supp\l=
\{g\in G\mid n_g\not= 0\}$.

Set
\bq
\wh\L_{\xi}=
\{\l\in
\wh{\wh \L}
\mid \forall C\in\ZZZ \quad \supp\l\cap\xi^{-1}([C,\infty[)
\mbox{ is finite }\}
\end{equation}

Then $\Lxi$ has a natural structure of  a ring.
We shall also need an analog of this ring with $\ZZZ$  replaced
by $\QQQ$.
That is, let
 $\wh{\wh \L}_\QQQ$
 be the set of all the functions $G\to\QQQ$
and
denote by
$\wh\L_{\xi,\QQQ}$
the subgroup of
$\wh{\wh \L}$
formed by all the $\l$ \sut~
for every $C\in\ZZZ$ the set
$\xi^{-1}([C,\infty[)\cap\supp\l$
is finite.  Then
$\wh\L_{\xi,\QQQ}$
is a ring, and $\Lxi$ is a subring.

A basic example of a Novikov ring is provided by
$G=\ZZZ$ and the standard inclusion
$\xi:\ZZZ\rInto\RRR$.
In this case $\Lxi$
is the ring of Laurent power series
with integral coefficients and finite negative part.

One of the main algebraic instruments of the paper are the groups
$\klxi$
and
$\klxiQ$.
We shall need also the analogs of Whitehead groups.
To introduce these groups
define first
a subset
$\pm G$
of $\L=\ZZZ G$
by:
$\pm G=\{\pm g | g\in G\}$.
Then
$\pm G\sbs \Lxi^\bu$
and this inclusion induces
a homomorphism
$\pm G\to \klxi$.
The image of $\pm G$
via this homomorphism
is denoted by $\pm \wh G$
and the group
$\klxi /\pm \wh G$
is denoted by $\wh {{\rm Wh}} (G, \xi)$.
Similarly the quotient
of
$\klxiQ$
by the image of
$\pm G$
is denoted by
$\wh {{\rm Wh}}_\QQQ (G, \xi)$.

\subsubsection{The Novikov complex and its simple homotopy type}
\lb{susu:novcom}

Let
$M$ be a closed connected \ma, and let $f:M\to S^1$ be a Morse map.
We assume that $f_*:H_1(M)\to H_1(S^1)=\ZZZ$
is epimorphic.
 Let $\PP:\wi M\to M$ be the universal
covering; then
 $f\circ \PP$ is homotopic to zero. The structure group
of this covering
(isomorphic to $\pi_1(M)$ ) will be denoted by $G$.
There is
a \ho~$\xi: G\to\ZZZ$, induced by $f$.
Set $\L=\ZZZ G$, then the corresponding
Novikov ring $\Lxi$ is defined.
Let $v$ be an $f$-gradient satisfying \TA~ (that is,
for every pair of zeros
$p,q$ of $v$
we have:
$W^{st}(p)\pitchfork W^{un}(q)$
where $W^{st}(p)$
and $W^{un}(q)$
are respectively stable and unstable
manifolds of $p$, resp. $q$).
Choose for every critical point $p$ of $f$ an orientation of the
stable \ma~of $p$
and a lifting of $p$ to $\wi M$. To this data one associates a
$\Lxi$-complex $ C_*(v)$,
called {\it the Novikov complex}
\sut~the number of free
generators of $C_k(v)$
equals  the number of critical points of $f$ of index $k$.
This chain complex has the same homotopy type as the complex
$$\cmxi= C_*^\D(\wi M)\tens{\L}\Lxi $$
where
$C_*^\D(\wi M)$
is the simplicial chain complex of $\wi M$
(see \cite{patou}).
Both
$C_*(v)$
and
$\cmxi$
are free based chain complexes
so a chain homotopy equivalence between them
has a naturally defined torsion in
$  \ove{K_1}(\Lxi)  $.
The bases in these complexes are determined by the choices of
orientations and
liftings of critical points
(resp. simplices)
to $\wi M$.
Changing these choices leads to the  changing of the torsion by adding 
an element
of
$\Im(G\to  \ove{K_1}(\Lxi)  )$.
Thus the torsion is really well defined in the group
$\whgxi$
(introduced above in the section \ref{susu:algnot}).

\subsubsection{Vector fields and closed orbits: terminology}
\label{susu:vecterm}

Let $w$ be a $\smo$ vector field on a closed manifold $M$.
The set of all closed orbits of $w$ is denoted by
$Cl(w)$.
A vector field $w$ will be called {\it Kupka-Smale } if
every zero and every closed orbit is hyperbolic, and the
stable manifold of every
zero of $v$
is transversal to the unstable
\ma~of every
zero of $v$.
For a Kupka-Smale vector field $w$ and a closed orbit
$\g\in Cl(w)$
let
$\ve(\g)\in \{-1,1\}$
denote
the index of the corresponding Poincare map; let
$m(\g)$
denote the multiplicity of $\g$.
Let $f:M\to S^1$ be a Morse map.
Gradient-like vector fields for $f$
(in the sense of the definition in \cite{milnhcob}, \S 3.1)
 will be also
called {\it $f$-gradients }.
The set of all Kupka-Smale $f$-gradients 
will be denoted by $\GG(f)$.

\subsubsection{Witt vectors and their logarithms}
\label{susu:wittlog}

Let $A$ be a ring with unit
 (non-commutative in general).
We assume that $\QQQ\sbs A$.
Let $\rho:A\to A$
be an automorphism.
Consider the
{\it $\rho$-twisted polynomial ring}
$A_\rho[t]$
(where the multiplication satisfies $at=t\rho(a)$
for $a\in A$)
and embed it into the corresponding
formal power series ring $P=A_\rho[[t]]$,
which is in turn a subring of the ring
$R=A_\rho((t))=
A_\rho[[t]][t^{-1}]$
of Laurent powers series with finite negative part.

The aim   of the Section \ref{s:wittvec}
is to construct a \ho~
$\LL$
defined on $K_1(P)$
and with values in an abelian group which we shall now describe.
For $n\in\NNN$
set $P_n=At^n\sbs P$.
Let $P_n'$
be the abelian subgroup of
$P_n$ generated by all the commutators $[x,y]=xy-yx$
where
$x\in P_k, y\in P_s, k+s=n$.
Set
\bq
\bar P_n =P_n/P_n', \quad P'= \prod_{n\geq 0} P_n', \quad \bar P
=
\prod_{n\geq 0} \bar P_n=
P/P'
\end{equation}
(In general $P'$
is not an ideal of $P$,
so $\bar P$
has only the structure of a $\QQQ$-vector space.)

The  homomorphism
$\LL$
 takes its values in
$\bar P$
and
 has in particular the following
logarithm-like
property:
\pa
{\it
For $x\in P$
the image of $tx$ in $\bar P$ is equal to $\LL([\exp tx])$.}
\pa
The main technical tool in the construction of $\LL$
is the group of Witt vectors in $P$.
A {\it Witt vector}
is an element of $P=A_\rho [[t]]$
of the form
$1+a_1t+a_2t^2+\dots$, with $ a_i\in A$.
The Witt vectors form a multiplicative subgroup $W$
of the full group of units $P^\bu$.
The image of $W$ in $K_1(P)$
will be denoted by $\wi W$. We first construct $\LL$
on the group $\wi W$ and then extend it to the whole of
$K_1(P)$ (section \ref{s:wittvec}).
The resulting extension 
$\LL:K_1(P)\to \bar P$
satisfies $\LL(K_1(A))=0$.

For our topological applications we need rather
a homomorphism
defined on the group $K_1(R)$.
The construction of such \ho~
follows immediately from the main theorem of the paper
\cite{ranpaj}
of A.Ranicki and the author.
The corollary which we need of the main theorem of \cite{ranpaj}
says:

\bepr[\cite{ranpaj}, Corollary 0.1]\label{p:cit_ranpaj}

The homomorphism
$\wi W\to \karzz$
induced by the inclusion
$W\rInto (\arz)^\bu$
admits a left inverse
$\wh B_2:\karzz\to \wi W$, \sut~
\been\item 
$\wh B_2$ vanishes on the image of $K_1(A)$ in $\karzz$
\item $\wh B_2$ vanishes on the image in $\karzz$
of the $1\times 1$-matrix
$(t)$.
\enen 
$\qs$
\enpr

Set $\gL=\LL\circ\wh B_2$; then the \ho~
$\gL:\karzz\to \bar P$
vanishes on the image of $K_1(A)$ in $\karzz$
and
the composition
$\karz\to\karzz\rTo^{\gL}\bar P$
equals to $\LL$.

\subsubsection{Logarithms in the context of group rings}
\label{susu:term_group}

The results about Witt vectors in power series rings
or Laurent series rings
will be applied to the Novikov completions of group rings.

Let  $G$ be a group and $\xi:G\to\ZZZ$ be an epimorphism.
Let $H=\Ker\xi$. Choose any
 element $t$
with $\xi(t)=-1$.
The ring  $\QQQ G$ is then identified naturally with
$\rho$-twisted Laurent polynomial ring
$A_\rho[t, t^{-1}]$, where $A=\QQQ H$,
and $\rho(x)=t^{-1}xt$
(so that the multiplication in this ring satisfies
$at=t\rho(a)$ for $a\in A$).
The Novikov ring
$\wh \L_{\xi, \QQQ}$
is identified with the twisted Laurent  series ring
$R=A_\rho((t)) $.

Let $\G$ be the set of conjugacy classes of $G$, and
let $\G_n$ denote the set of all the conjugacy classes
contained in $G_{(n)}=\xi^{-1}(n)$.
Consider the $\QQQ$-vector space
$\QQQ \G_n $
spanned by $\G_n$, and set
$\gG   =
\prod_{n\leq 0} \QQQ
\G_n$.

Lemma
\ref{l:pbar_gamma}
 shows that $\gG\approx \bar P$
where $P=(\QQQ H)_\rho [[t]]$.
Thus we obtain a homomorphism
$\klxiQ \to \gG$.
Up to natural  identifications
it is still the same \ho~
$\gL: K_1(R)\to \bar P$,
so we keep the same symbol $\gL$
to denote it.
(The advantage of the notation
$\gL:\klxiQ \to \gG$
is that the objects
$\klxiQ, \gG$
do not depend on the particular choice of the element  $t$.)

In subsection \ref{su:grouprings}
we show that $\gL$
factors through
$\whgxiQ=
\klxiQ/\pm\wh G$.
Composing with the embedding $\Lxi\sbs \LxiQ$
we obtain now a \ho~
from
$\whgxi$
to $\gG$.
We shall denote it by the same symbol
$\gL$
since \noconf.

\subsubsection{Non abelian eta function for gradient flows}
\lb{susu:nonabeta}

Let $f:M\to S^1$
be a Morse map and let $v$ be an $f$-gradient.
Let $\g$ be a closed
orbit of $(-v)$. There is no natural choice
of a base point in $\g$,
and the class of $\g$ in $\pi_1(M)$
is ill-defined, but the class of $\g$
is well-defined as the element of the set $\G$
of conjugacy classes
of $\pi_1(M)$.
This class
will be denoted by $\{\g\}$.
We are about to introduce the eta-function
counting the classes $\{\g\}$ corresponding to
 the closed orbits of $-v$.

Let $\xi=f_*:\pi_1(M)\to\ZZZ$.
We assume that $\xi$ is epimorphic.
 Applying the terminology of the subsection
(\ref{susu:term_group})
 we obtain
the subsets
$\G_n\sbs \G$
and the rational vector space
$\gG$.
Assume now that $v$ is a Kupka-Smale vector field.
Here is the definition of  the non-abelian eta function of $-v$:
\begin{equation}
\eta_L   (-v)=\sum_{\g\in Cl(-v)} \frac {\ve(\g)}{m(\g)} \{\g\}
\in \gG
\lbl{f:etadef}
\end{equation}

(Using the Kupka-Smale condition, it is not difficult to check that
for a given $n$ there is only a finite number
of closed orbits $\g$ with
$\{\g\}\in
\G_n $, and, therefore, the right hand side of
(\ref{f:etadef})
is indeed a well-defined element of $\gG$.)

 Recall the notation
$
\cmxi=
C_*^\D(\wi M)\tens{\L}\Lxi $
from subsection \ref{susu:novcom}.

\pa
{\bf Main Theorem }

\it There is a subset
$\GG_0(f)\sbs \GG(f)$
with the following properties:
\been
\item
$\GG_0(f)$ is open and dense in $ \GG(f)$
\wrt~$C^0$ topology
\item
For every
$v\in \GG_0(f)$
there is a chain \hot~equivalence
\begin{equation}
\phi: C_*(v)\rTo^
\sim
\cmxi
    \lbl{f:ntzeta1}
\end{equation}
 \sut~
\begin{equation}
\gL(\tau(\phi ))=-\eta_L   (-v)    \lb{f:ntzeta2}.
\end{equation}
\end{enumerate}

\pa

\rm

\subsubsection{Main lines  of the proof}
\lb{su:ideaprf}

In the  paper
\cite{pastpet},
we constructed for a given Morse map $f:M\to S^1$
a special class of $f$-gradients which form a subset,
$C^0$-open-and-dense
in the subset of all $f$-gradients satisfying \TA.
For any $f$-gradient $v$ in this class the boundary  operators in the 
Novikov
complex, associated to $v$, are not merely power series, but rational
 functions.
The main construction in the proof is a group \ho~
$h(v)$
between \hog~groups of certain spaces,
 (see  \S 4 of \cite{pastpet})
which we call
   {\it homological gradient descent}.

To explain the roots  of this notion,
let $\l\in S^1$
be a regular
value of $f$, and consider the integral
curve $\g$ of $(-v)$ starting at
$x\in V=f^{-1}(\l)$.
If this curve does not converge to a critical point of $f$, then
it will intersect again the \sma~$V$ at some point, say, $H(x)$.
The map $H$ is thus a (not everywhere defined ) smooth map
of $V$ to $V$, and the {\it homological gradient descent operator}
is a substitute for the \ho~induced by $H$ in \hog.
This operator is defined for a $C^0$-generic $v$-gradient satisfying 
\TA.

Note that the closed orbits of $v$ correspond to
 the periodic points of the map
$H$. This makes possible to give a formula
for $\eta_L$  in terms
of the homological
gradient descent operator
(see the formula (\ref{f:htau})).
On the other hand it turns out that there is a  homotopy equivalence
(\ref{f:ntzeta1})
such that its torsion is computable
in terms of homological gradient descent operator
(see (\ref{f:tors})).
Comparing these two expressions we obtain the theorem.

\subsection{Acknowledgements}
\label{su:ackno}

Many thanks are due to the referee
for careful reading of the manuscript, and for
pointing out two   errors in the first version of the
paper.

I am happy to thank  A.Ranicki
for the attention he paid to this work
and for our collaboration on the subject of the Whitehead groups
of the Novikov rings.
Our joint work 
 has lead to the paper \cite{ranpaj}, the results of which
allowed a strengthening of
 the results of the first version of the present paper.

I am grateful to  V.Turaev
for valuable discussions and
to A.Fel'shtyn for useful information 
 about zeta functions.

After this paper has been accepted for publication
 I received a very interesting preprint of D.Sch\"utz
\cite{schutz}.
Sch\"utz's  work  contains in particular a generalization
of the results of the present paper to the case of 
closed 1-forms within an irrational cohomology class.
Many thanks to D.Sch\"utz for having sent me his paper.

A part of the work was done  during
my stay at IHES in December 1998
and May 1999. I am grateful to IHES for hospitality and financial 
support.

\section{Witt vectors and homomorphism $\LL$}
\label{s:wittvec}

We shall be working here with the notation of
\ref{susu:wittlog}.
Thus $A$ is a ring containing $\QQQ$, $\rho:A\to A$ is
an automorphism, $P=\arz$
is the  corresponding twisted power series ring,
$R=\arzz$
is the  twisted Laurent  series ring, etc.
For $x\in R,   x=\sum_ix_it^i$
set $\nu(x)=\min\{k\in \ZZZ \mid x_k\not= 0\}$.
Set $P_+=\{x\in P|\nu(x)>0\}$.
There is a natural topology in $P$, namely
$a_n\to a$ if $\nu(a_n-a)\to \infty$.
Note that the closure $\ove{[P,P]}$
is in $P'$.
Our aim in this section is to construct the
homomorphism $\LL:K_1(P)\to \bar P$.
 We do it in three steps:
first we define a logarithm-like homomorphism
on the group of Witt vectors $W$, then
we factor it through
$\wi W=\Im(W\to K_1(P))$ ;
finally we extend it to the whole of $K_1(P)$
using the computation of
$K_1(P)$ from \cite{patou}.
(Note by the way, that we do not use here the results
of \cite{ranpaj}.)

\subsection{The
Baker-Campbell-Hausdorff-Dynkin formula in $A[[t]]$}
\lb{su:chd}

Define the exponent and logarithm
 as usual by the following formulas

\begin{gather}
\log (1+\mu)=\mu-\frac {\mu^2}2+\frac {\mu^3}3-\dots+(-1)^{n-1}\frac 
{\mu^n}n
+\dots , \mbox{ where }   \mu\in P_+,\lbl{loga} \\
\exp y=e^y=1+y+\frac {y^2}2+\dots
\frac {y^n}{n!}+\dots \mbox{ where }    y\in P_+ \lbl{expo}
\end{gather}

It is clear  that  both  series
converge in $P$.
A standard computation shows that
$\log\exp y = y$ for $y\in P_+$,
and $\exp \log x =x$ for $x\in W$.
For $X,Y\in P_+$ set

\begin{equation}
\zeta(X,Y)=
\log (e^Xe^Y)\in \arz\lbl{zeta}
\end{equation}

The well-known
Baker-Campbell-Hausdorff-Dynkin computation
(see for example
\cite{postnlie}, Lecture 6, or
\cite{bour}
Ch. 6, \S 4 )
implies:
\begin{equation}
\z(X,Y)=
\sum_{n=1}^\infty
\z_n(X,Y)      \lbl{f:dzeta}
\end{equation}
where
$\z_n(X,Y)$
is a Lie polynomial in $X, Y$ of degree $n$.
The power series in the right hand side of
(\ref{f:dzeta})
converges since
$\nu(\zeta_n(X,Y))\geq n$. Moreover, for $n>1$ we have:
$\zeta_n(X,Y)\in [P,P]\sbs P'$.
Therefore if $\mu,\nu$
are the units in $W$ of the form
$e^X, e^Y$
with
$X,Y\in P_+$, we have
$\log(\m\n)-\log \m -\log \n \in P'$.
Since every element of $W$ is of the form
$e^X$
with $X\in P_+$, we obtain that the composition
\begin{equation}
W\rTo^{\log} \arz =P \rTo^{pr} \bar P
\end{equation}
is a group \ho, which will be denoted by $\Log$.
In the next section we shall see that this \ho~
can be factored through the group
$\wi W =\Im (W\to K_1(P))$.
The first step in the proof is the next lemma.
(Note, that the symbol $[P^\bu, P^\bu]$
denotes  the commutator of the group $P^\bu$
in the group-theoretic sense, that is
$[P^\bu, P^\bu]$
is the subgroup of $P^\bu$
generated by the set
$
\{ ghg^{-1}h^{-1} | g,h \in P^\bu\}$.)

\bele\lb{l:comm}
The \ho~ $\Log$
vanishes on the subgroup
$W\cap [P^\bu, P^\bu]$
of $W$.
\enle

\Prf
Let $m\in W$.
The elements of the form $\a m\a^{-1}$,
where
$\a\in A^\bu$
will be called
 {\it $A$-conjugate} to $m$.
It is obvious that the value of
$\Log $
on $A$-conjugate elements is the same.

Let $X\in W\cap [P^\bu, P^\bu]$.
Assume first that $X=
\kom ab$
where
$a,b\in P^\bu   $. Write $a=\a m, b=\b n$
where $\a,\b\in A^\bu$
and $m,n\in W$.
A direct computation gives:
$X=
\kommm {m}{n}\cdot \kom {\a}{\b}  $,
where
$m',m''$
are $A$-conjugate to $m$, and
$n',n''$
are $A$-conjugate to $n$.
Since $X\in W$ we have $\kom {\a}{\b}  =1$, and
$\Log (X)=0 $.
The general case ($X$ is a product of several commutators)
is done by an easy induction
argument.
$\qs$

\subsection{Homomorphism $\LL$}
\lb{su:lnma}
Apply the results of the previous section to the ring
$\AA_n=\MM_n(A)$  of $n\times n$-matrices over $A$.
We have the  power series ring
$\PP_n=\AA_n[[t]]$, the
multiplicative subgroup
$\UU_n=
1+t\PP_n\sbs
\PP_n^\bu
 $
of the group of units,
the  Laurent power series ring
$\RR_n$ etc.
By the results of the preceding subsection
~the composition
\begin{equation}
\UU_n\rTo^{\log}\AA_n[[t]]
=\PP_n
\rTo{pr}\bar\PP_n
\end{equation}
is a group \ho.
We have the trace map
$\Tr:\PP_n\to A[[t]]$
defined by
$\Tr(\sum_i\a_i t^i)=
\sum_i \Tr(\a_i) t^i$.
It is easy to check that $\Tr$ factors through the quotients
to define correctly the group \ho
~$\bar\PP_n\to \bar P$,
which we shall denote by the same symbol $\Tr$.
The composition $\Tr\circ pr\circ   \log$
will be denoted by
$\LL_n:\UU_n\to \bar P$;
it is a group \ho.
It is easy to check that the map
$\LL_n$
behaves well \wrt~
the
natural stabilization map
$\UU_n\rInto \UU_{n+1}$.
Set
$\UU=\limdir ~\UU_n$
and obtain
a \ho
~$\UU\to \bar P$,
which will be denoted by
$\LL$.
Let $\wi\UU$ be
 the image of $\UU$
via the \ho
~$\UU\rInto GL(P)\to  K_1(P)$.
It follows from the Lemma \ref{l:comm}
that $\LL$ factors through $\wi\UU$;
we keep the same notation
$\LL$ for the resulting \ho~
$\wi\UU\to \bar P$.
It is obvious that
$\wi\UU\supset \wi W=\Im(W\to K_1(P))$
 and
that the restriction of $\LL$ to $\wi W$
coincides with $\Log$.

\bepr\label{p:uandw}
\been\item
$\wi\UU=\wi W$.
\item
The inclusions
$W\rInto^i P,~ A\rInto^j P$ and the projection
$P\rTo^p A$
induce the following  split exact sequence \bq
0\rTo\wi W\rTo^{i_*} K_1(P)
\pile{\rTo^{p_*}
\\
\lTo_{j_*}}
K_1(A)\rTo 0 \lbl{f:splitpower}
\end{equation}
\enen
\enpr
{\it Proof. }\quad
A slightly weaker version of this proposition
is contained in
\cite{patou}, Lemma 1.1; it is based on the
argument  due essentially to A.Suslin.
Let us 
consider an element
$Z\in K_1(P)$ with $p_*(Z)=0$.
Then $Z$ can be reduced to a square matrix $S$
of the form
$1+\mu$
where $\mu$
has the coefficients in $P_+$.
Each diagonal element of $S$
is then invertible
in $P$ and using the standard Gauss elimination method
we can reduce $S$
to a diagonal matrix $\D$
without changing the class of $S$ in $K_1$.
The class of $\D$
in $K_1$
belongs obviously to $\wi W$.
$\qs$

In particular, $\LL$
can be extended in a functorial way
to a \ho~
of the whole of
$K_1(P)$  to $ \bar P$:
just define $\LL$
to be zero on the image of $j_*$.

Note that $\Im\LL=\bar P_+$.

\section{Laurent series rings and homomorphism $\gL$}
\label{s:laurent_log}

As we have already mentioned (see \ref{susu:wittlog})
the homomorphism $\LL$ gives rise to the \ho~
$\gL:K_1(R)\to \bar P$.
We shall work mainly with the case of group rings. 
This case has some particularities which we shall now explain.

\subsection{ The case of group rings}
\label{su:grouprings}

Here we consider a group $G$ and an epimorphism
$\xi:G\to\ZZZ$.
We shall work here with the terminology of
\ref{susu:wittlog}. Thus
$H=\Ker\xi, \quad \G$
stands for the
set of conjugacy classes of $G$,
 $\G_n$ is the set of all the conjugacy classes
contained in $G_{(n)}=\xi^{-1}(n)$.
Recall also the abelian group
$\gG=
\prod_{n\leq 0} \QQQ
\G_n$.
Set 
$\gG_+=
\prod_{n < 0} \QQQ
\G_n$.
 We choose and fix any element $t\in G_{(-1)}$.

\bele\label{l:pbar_gamma}
There is an  isomorphism $\bar P_+\rTo^i \gG_+$, \sut~the
 following diagram
is commutative
\begin{diagram}
\bar P_+ & \rTo^i & \gG_+ \\
 &\luTo(1,2)   \ruTo(1,2)&   \\
    & P &  \\
\end{diagram}
Here the oblique arrows are natural projections
( they are only abelian groups homomorphisms, not ring homomorphisms)
\enle
\Prf
It is convenient to use another description of $\gG$,
made by analogy with
$\bar P$.
Namely, let $n\geq 0$.
Let $\bar R_n'$
be the abelian subgroup of $R_n$,
generated by
all the commutators $[x,y]$
where
$x\in R_k, y\in R_s, k+s=n$
(we do not assume that $k, s$ are positive).
Set $\bar R_n=R_n/R'_n$
and consider an abelian group ${\bar R}=
\prod_{i=0}^\infty \bar R_i$,
so that
$\bar R$
is a quotient of $R$. Then $\bar R$
is identified in an obvious way with
$\gG$.
There is an inclusion
$P_n'\sbs R_n'$, and the natural projection
$\pi_n:\bar P_n\to \bar R_n$,
and to prove the lemma  it suffices
to prove that $\pi_n$
is an isomorphism for $n>0$.
The automorphism
$\rho:A\to A$
extends to  an automorphism
$\rho:P\to P$
of the ring $P$
by the following formula:
$\rho(\sum_ia_it^i)=\sum\rho(a_i)t^i$.
The set $P'$
is preserved by $\rho$, so $\rho$
determines an isomorphism of abelian groups
$\bar \rho:\bar P\to \bar P$.
Note that for $n>0$
the homomorphism
$\bar \rho_n:\bar P_n\to \bar P_n$
equals $\id$.
(Indeed, let
$a\in P, a=bt$ with
$b\in P_{n-1}$.
Write
$a-\rho(a)=bt-\rho(bt)=t\rho(b)-\rho(bt)=
\rho(tb-bt)$.
The last element is obviously in $\bar P$.)
Now we can prove the proposition.
Let $xy-yx\in R_n'$,
where $x=at^p, y=bt^{-q}$
with
$a,b\in A_0, p,q\geq 0$
and
$p-q=n>0$.
Consider the element
$yx=t^{-p}(t^pbt^{-q}a)t^{p}
\in P_n$.
Its class in
$\bar P_n$
equals the class of
$t^pbt^{-q}a$, which equals the class of
$at^pbt^{-q}=xy$
and we are over. $\qs$

Thus the homomorphism
$\gL: K_1(R)\to \bar P$
can be identified with a homomorphism
$K_1(\LxiQ)\to \gG$,
which will be denoted by the same symbol $\gL$.
It is easy to see that $\gL$
does not depend on a particular choice of $t\in G_{(-1)}$
so it depends only on
$G$ and $\xi$.

It follows now from Proposition 
\ref{p:cit_ranpaj}
that $\gL$ vanishes on the subgroup
$\pm \wh G$ of $K_1(\LxiQ)$.
Composing with the embedding 
$\Lxi\sbs \LxiQ$
we obtain thus a homomorphism
$\gL:\klxi\to\gG$,
which factors through
$\whgxi$
and the image of which is contained in $\gG_+$.

\section{ Application of the logarithms
to the study of the groups $K_1(P), K_1(R)$}
\label{s:applic_log}

We have seen that the group $\klxi$
admits a natural splitting.
The last two direct summands of this splitting,
namely
the Whitehead group
$K_1(\ZZZ H)$
and the $Nil$-group
of $H$
are well known in $K$-theory.
This is not the case with the group
$\wi W$, the structure of which is not yet well understood.
The situation is not clear even 
if we consider $K_1$
of the rational Novikov ring
$\LxiQ$.
The logarithm maps, defined above, can shed some light on
 the structure of these groups. 
In the rest of this section
we shall be working in the terminology of
\ref{susu:wittlog}.
Thus
$\QQQ\sbs A, P=\arz, R=\arzz$.
We shall describe explicitly the kernel
of the map $\Log:W\to \bar P$. This provides a family of
 non-vanishing elements
in the group $\wi W=\Im (I:W\to K_1(P))$,
since $\Log=\LL\circ I$ and
the elements of $W$,  which are not in $\Ker \Log$
are not in the $\Ker I$.
(Note that since the map
$\wi W\to K_1(R)$
is a split monomorphism,
we get automatically the information about the group
$\Im(W\to K_1(R))$.)

\subsection{Power series}
\label{su:log_pow}

We have already mentioned that $K_1(P)\approx
\wi W\oplus K_1(A)$.
Thus the computation of
$K_1(P)$
is reduced to the computation
of $K_1(A)$ and the image of $W$
in $K_1(P)$.
Recall from Section 1  a commutative diagram:
\begin{diagram}
W   &  \rTo^I  &  K_1(P) \\
&  \rdTo_{\Log}  & \dTo^\LL  \\
&  & \bar P_+ \\
\end{diagram}
Thus $\Ker I\sbs \Ker \Log$.
In particular 
$W_1=W\cap [P^\bu, P^\bu]$
is in 
$\Ker \Log$.
\footnote{
Recall, that the symbol $[P^\bu, P^\bu]$
denotes  the commutator of the group $P^\bu$
in the group-theoretic sense, that is
$[P^\bu, P^\bu]$
is the subgroup in $P^\bu$ generated by the set
$\{ ghg^{-1}h^{-1} | g,h \in P^\bu\}$.
}
Since $\Log$
is continuous in the natural topology of $P$,
the closure
$\ove{W_1}$
is also in $\Ker\Log$.
It turns out that for the rings $A$ satisfying the property
(*) below
the kernel of $\Log$
coincides with
$\ove{W_1}$.
\pa
(*) \qquad Every element of $A$
is a sum of invertible elements
\pa
(This condition is satisfied, for example, for
the rational group rings.)
Let us denote by $\ove{\Log}$
the homomorphism
$W/\ove{W_1} \to \bar P_+ $
derived from $\Log$.

\beth
\label{t:kerl}
If $A$ satisfies
(*), then
$\ove{\Log}: W/\ove{W_1} \to \bar P_+ $
is an isomorphism.
\enth
\Prf
Our task is  to check the injectivity.
Let $\a\in W$ belong to $ \Ker\ove{\Log}$.
We shall find a sequence of elements
$\b_n\in W_1$
such that
$\nu(1-\a \cdot \b_1\cdot \dots \cdot \b_n)\to\infty$
and
$\nu(1-\b_n)\to\infty$
as
$n\to\infty$.
This will prove our result since the element
$\b= \lim_{n\to\infty} \b_1\cdot \dots \cdot \b_n$
is in $\ove{W_1}$
and
$\a=\b^{-1}$.
The existence of such a sequence of elements follows from the next 
lemma.

\bele\lb{l:induk}
Let
$\a\in W, \nu(1-\a)=n$.
Assume that
$\nu(\Log(\a))\geq n+1$.
Then there is an element
$X\in W_1, \nu(1-X)=n$,\sut~
$\nu(1-\a X)\geq n+1$.
\enle
\Prf
Let $\a=1+x+ O(t^{n+1}) $
(where $x\in P_n$ and $O(t^{n+1}) $ stands for
"terms of degree $\geq n+1$").
Then
$\log\a=
x+ O(t^{n+1}) $, and
thus $x\in P_n'$.
Let us denote by $P_n^\bu$
the subgroup of $P^\bu$
consisting of all elements $\xi=\xi_0+\xi_1+\dots+\xi_n+\dots$, 
\sut~$\xi_0=1$ and $\xi_i=0$
for
$0<i<n$.
Our lemma will be proved if we show the following:

\bele
Let $z\in P_n', n>0$.
There is
$\xi\in P_n^\bu\cap W_1$, \sut~$\xi_n=z$.
\enle
\Prf
Let us call an element $z\in P_n$
{\it regular}
if $z=\xi_n$ for some
$\xi\in P^\bu_n\cap W_1$.
Note that a sum of regular elements is again regular.
Thus it suffices to prove that for any 
elements $x\in P_s, y\in P_r$ with $s+r=n$ and
$s,r\geq 0$ the 
commutator
$xy-yx$
is regular.
We consider two cases:
\been\item
One of $\nu(x), \nu(y)$ is zero.
Assume for example that $\nu(x)=0$.
In view of the condition $(*)$
we can assume that $x$ is invertible.
Then
$xy-yx=\xi_n$
where
$\xi=
x(1+yx)x^{-1}(1+yx)^{-1}$.
\item
Both $\nu(x), \nu(y)$ are strictly positive.
Let 
$\xi=
(1-x)(1-y)(1-x)^{-1}(1-y)^{-1}$.
Expanding this expression, one shows easily that the first non-zero
term of $\xi-1$
equals $xy-yx$.
\enen $\qs$

\begin{coro}\label{c:laur_ser}
If $A$ satisfies $(*)$ then 
$\ove{\Ker(W\to K_1(R))}=
\ove{\Ker I}=
\ove{W_1}= \Ker \Log$. $\qs$
\end{coro}

\section{Counting closed orbits: preliminaries}
\label{s:cloprel}

The proof of Main Theorem
uses the techniques from
the author's papers \cite{pamrl}, \cite{pastpet}, \cite{pafest}.
In this section we recall the
necessary material following essentially \cite{pafest}.
The main objects in this section are 
{\it $C^0$-generic gradients} of  a Morse map $f:M\to S^1$.
Subsections \ref{su:mofgr} and \ref{su:dthd}
contain recollections of some
more or less
standard material
in Morse theory,
centered around Morse complexes and handle decompositions.
The notion of {\it gradient descent homomorphism}
associated to a given $f$-gradient
is the contents of the section
\ref{su:celappr}.
This homomorphism is
defined only for some special class of gradients.
The $C^0$-generic condition $(\gC)$
on  vector fields
which suffices to define such a homomorphism
is the contents of the section \ref{su:cc}.
Recall from \cite{patou}
that for every $f$-gradient $v$ satisfying
\TA~
there is a chain homotopy equivalence
$\wi C_*(v)\to
\cmxi$.
In the paper \cite{pafest}
we proved that for a $C^0$-generic gradient
satisfying
\TA~ the chain homotopy equivalence above
can be chosen so that its torsion
is computable in terms of the  gradient descent homomorphism
associated to $v$.
The resulting formula 
(\ref{f:tors})
for the torsion is given
in the subsection
\ref{su:circle}.
The proof of Main Theorem will be obtained in the Section 
\ref{s:prfclo}
 by
 identifying the right hand side of
(\ref{f:tors})
with the non-abelian 
eta function of the gradient flow.

\subsection{Morse functions and their gradients}
\lb{su:mofgr}

Let $v$ be a $C^1$ vector field on a \ma~ $M$.
The value of the integral curve of $v$ passing by $x$ at $t=0$
will be denoted by $\g(x,t;v)$.
The set of all critical points of a Morse
 function $f$ will be denoted by $S(f)$; the set
 of all critical points of $f$ of index $k$
will be denoted by $S_k(f)$.

Let $f:W\to [a,b]$  be a Morse function on a cobordism
$W$,
so that $f^{-1}(b)=\pr_1 W, f^{-1}(a)=\pr_0 W,
\pr W=\pr_1 W\sqcup \pr_0 W$.
We shall call gradient-like vector fields of $f$
(in the sense of the definition in \cite{milnhcob}, \S 3.1)
also {\it $f$-gradients}.
The set of all $f$-gradients will be denoted by $G(f)$.
Let $v\in G(f)$.
Set $U_1=\{x\in \pr_1 W|\g(x,\cdot;v) \mbox{ reaches } \pr_0 W\}$.
Then $U_1$ is an open subset of $\pr_1 W $
and the gradient descent along the trajectories of $v$
determines a diffeomorphism of the  open subset
$U_1$ of $\pr_1 W$
onto an open subset $U_0$
of $\pr_0 W$.
This diffeomorphism will be denoted by
$\stv$
and we abbreviate
$\stv (X\cap U_1)$ to
$\stv (X)$.
For $x\in S(f)$
we denote by $D(x,v)$
the
{\it descending disc of $x$},
that is the set of all $y\in W$, \sut
~$\g(y,t;v)\rTo_{t\to\infty} x$.

We say that
$v$ satisfies {\it Almost Transversality Condition},
if
   $$\big( x,y \in S(f) ~\&~
\ind x \leq \ind y \big) \Rightarrow \big(
D(x,v) \pitchfork D(y,-v) \big)$$
\vskip0.1in
We say that
$v$ satisfies {\it Transversality Condition},
if
   $$\big( x,y \in S(f) \big)
\Rightarrow \big(
D(x,v) \pitchfork D(y,-v) \big)$$
\pa
The set of all $f$-gradients satisfying
\TA, resp. \ATA~
will be denoted by $GT(f)$, resp. by $GA(f)$.

\pa

A very useful class of Morse functions is that of {\it ordered 
functions},
(which is a slightly wider class than the class of 
self-indexing functions of
S.Smale).
We say that $f:W\to[a,b]$ is {\it ordered} Morse function
with an {\it ordering sequence} $(a_0,\dots,a_{n+1})$, if
$a=a_0<a_1<\dots<a_{n+1}=b$ are regular values of $f$
\sut~ $S_i(f)\subset
 f^{-1}(]a_i,a_{i+1}[)$.
(Here $n=\dim W$.)

In many cases the considerations involving
 gradients of arbitrary Morse functions
can be reduced to gradients of ordered Morse
functions. (This is for example the case with the definition of the
Morse-Thom-Smale-Witten complex of a Morse function.)

 A Morse function $\phi:W\to[\a,\beta]$
is called {\it adjusted to $(f,v)$}, if:

 1) $S(\phi)=S(f)$, and $v$ is also a $\phi$-gradient.

 2) The function $f-\phi$ is constant in a \nei~
of $\pr_0W$, in a \nei~ of $\pr_1 W$, and in a \nei~ of each
point of
$S(f)$.

One can show that for an arbitrary Morse function $f$
and an $f$-gradient satisfying the \ATA,
there is an ordered Morse function $g$, 
adjusted to $(f,v)$.

\pa

\subsection{ $\d$-thin handle decompositions}
\lb{su:dthd}

In this subsection $W$ is a riemannian cobordism
of dimension $n$, $\fcob$ is a Morse function on $W$, and $v$ is an
$f$-gradient.
We denote $W\sm \pr W$ by $W^\circ$.

Let $p\in W^\circ$. Let $\d>0$.
Assume that for some $\d_0>\d$ the restriction of
the exponential map
$\exp_p:T_pW\to W$ to the disc $B^n(0,\d_0)$
is a \dfm~on its image.
Denote by $B_\d(p)$ (resp. $D_\d(p)$) the riemannian open ball
(resp. closed ball)
of radius $\d$
centered in $p$. We shall use the notation
$B_\d(p), D_\d(p)$ only when the assumption above on $\d$ holds.
Set
\begin{gather*}
B_\delta(p,v)=\{x\in W~\vert~ \exists
t\geq 0 :
 \gamma (x,t;v)\in  B_\delta (p)\}\\
D_\delta(p,v)=\{x\in W~\vert ~\exists
t\geq 0 :
\gamma (x,t;v)\in  D_\delta (p)\}
\end{gather*}

We  denote by $D (\indl s;v)$
the union of   $
D(p,v)$ where $p$ ranges over critical
points of $f$ of index $\leq s$.
We denote by $B_\d (\indl s;v)$, resp.
by $D_\d(\indl s;v)$
the union of
$B_\d(p,v)$, resp. of
$D_\d(p,v)$, where $p$ ranges over critical
points of $f$ of index $\leq s$.
We shall use
similar
notation like $D_\d(\inde s;v)$ or
$B_\d(\indg s;v)$, which is now clear
without special definition. Set
$$C_\d(\indl s; v)=W\sm B_\d(\indl {n-s-1};-v)$$

Let $\phi:W\to[a,b]$ be an ordered Morse function with an
ordering sequence
$(a_0<a_1\dots <a_{n+1})$. Let $w$ be a $\phi$-gradient.
Denote $\phi^{-1}\([a_i,a_{i+1}])$ by $W_i$.

\begin{defi}
We say that $w$ is {\it $\d$-separated \wrt~
$\phi$ } (and the ordering sequence $(a_0,\dots,a_{n+1})$), if

i) for every $i$ and every $p\in S_i(f)$ we have
$D_\d(p)\sbs \Wkr_i$;

ii) for every $i$ and every $p\in S_i(f)$
there is a Morse function
$\psi:W_i\to[a_i,a_{i+1}]$, adjusted to
$(\phi\mid W_i, w)$
and a regular value $\l$ of $\psi$ \sut~
$$D_\d(p)\sbs\psi^{-1}(]a_i,\l[)$$
and for every $ q\in S_i(f), q\not= p$ we have
$$D_\d(q)\sbs \psi^{-1}(]\l,a_{i+1}[)$$

\end{defi}

We say that $w$ is {\it $\d$-separated }if
it is $\d$-separated \wrt~some ordered Morse
function $\phi:W\to [a,b]$, adjusted to $(f,v)$.
Each $f$-gradient satisfying \ATA~ is $\d$-separated
for
some $\d>0$.

\bepr [\cite{pafest}, Prop. 3.2, 4.1]
\lb{p:ddd}

If $v$ is  $\delta_0$-separated, then
$\forall
\delta\in[0,\delta_0]$ and
$\forall
s: 0\leq s\leq n$
\been
\item $D_\delta (\indl s;v)$         is compact.
\item ${\bigcap} _{\theta>0} B_\theta (\indl s;v)
= D (\indl s;v)$
\item
 $B_\delta (\indl s;v)
=\Int D_\d(\indl s; v)$
and
$ D_\d(\indl s; v)
\sbs
 C_\d(\indl s; v)$.
\item
$H_*(
 D_\d(\indl s; v)\cup \pr_0 W,
 D_\d(\indl {s-1}; v)\cup \pr_0 W)$
equals $0$ if $*\not= s$
and is a free module generated by the
classes of the descending discs $D(p,v)$
with $p\in S_s(f)$.
\end{enumerate}
\enpr

Thus the  collection of descending discs $D(p,v)$ form
a sort of stratified manifold, and
the open sets  $B_\d(v)$
 form  a family  of $\d$-thin
neighborhoods of this manifold.

We shall often denote
$D_\d(\indl s;v)$
by $W^{\{\leq s\}}$
if the values of $v,f,\d$
are clear from the context.

\subsection{Condition $(\gC)$}
\lb{su:cc}

In this subsection we recall the condition
$(\gC)$
on the gradient $v$.
If this condition holds, the gradient descent map, corresponding
to
$v$ can be endowed with a structure, resembling closely
the cellular maps between $CW$-complexes.
We shall
explain this cellular-like structure
in the section
\ref{su:celappr}.
We begin by stating the condition $(\gC)$.
Let $f:W\to[a,b]$ be a Morse function on a
riemannian
cobordism $W$,
$v$ be an $f$-gradient.

\bede\lb{d:cc}
(\cite{pafest}, Def. 4.5)

We say that
{\it $v$ satisfies condition $(\gC)$}
if there are objects 1) - 4), listed below,
with the properties (1 - 3) below.
\pa
{\it Objects:}
\pa
\begin{enumerate}
\item[1)] An ordered Morse function $\phi_1$
on $\pr_1 W$ and a $\phi_1$-gradient $u_1$.
\item[2)]  An ordered Morse function $\phi_0$
on $\pr_0 W$ and a $\phi_0$-gradient $u_0$.
\item[3)]  An ordered Morse function $\phi$
on $ W$
adjusted to $(f,v)$.
\item[4)]  A number $\d>0$.
\end{enumerate}
\pa
{\it Properties:}
\pa
(1)\quad  $u_0$ is $\d$-separated \wrt~ $\phi_0$,
$u_1$ is $\d$-separated \wrt~ $\phi_1$, $v$ is
$\d$-separated \wrt~ $\phi$.
\pa

 \begin{multline*}(2)\qquad    \stv \bigg(C_\d(\indl j;u_1)\bigg)
\cup
\bigg(D_\d\(\indl {j+1};v\)\cap \pr_0 W\bigg)
\sbs                \\
B_\d(\indl j,u_0)
\mbox{ for every } j
\end{multline*}

\begin{multline*}(3) \qquad     \st v\bigg(C_\d(\indl j;-u_0)     
\bigg)
\cup
\bigg( D_\d(\indl {j+1}; -v)\cap \pr_1 W\bigg)
\sbs                  \\
B_\d(\indl j;-u_1)
\mbox{ for every } j
\end{multline*}
$\qt$
\end{defi}

The set of all $f$ gradients satisfying $(\gC)$ will be denoted by
$GC(f)$.

\pa

\begin{theo}\lb{t:cc}
(\cite{pafest}, Th. 4.6)
$GC(f)$ is open and dense in $G(f)$ \wrt~  $C^0$ topology.
Moreover, if $v_0$ is any $f$-gradient then one can choose a $C^0$
 small perturbation $v$
of $v_0$ \sut~ $v\in GC(f)$ and
$v=v_0$ in a \nei~ of $\pr W$.
\end{theo}

\subs{Homological gradient descent and a cellular approximation for
$\stv$}
\label{su:celappr}

As we have already mentioned the application
$\stv$
is not everywhere defined. But if the gradient $v$
satisfies the condition
$(\gC)$,
we can associate to $v$
some family of continuous maps which plays the role of
"cellular approximation" of $\stv$,
and a homomorphism
$\HH(-v)$
({\it homological gradient descent})
which is a substitute for the
homomorphism induced by $\stv$ in homology.

So let $v$ be an $f$-gradient satisfying $(\gC)$.
It will be convenient to denote
$\pr_0 W$ by $V_a$, and $\pr_1 W$
by $V_b$.
We have the corresponding functions
$\phi_0, \phi_1$
on
$V_a$, resp.
$V_b$; and their
gradients $u_0, u_1$
give rise to the corresponding handle-like filtrations

\bq
V_a^{\{\leq s\}} = D_\d(\indl s; u_0),\quad
V_b^{\{\leq s\}} = D_\d(\indl  s; u_1) \lbl{f:vas}
\end{equation}

Consider the set $Q_s$
of all $x\in \vbs$
where
$\stv$ is {\it not} defined.
Equivalently,
$$
Q_s=\{x\in \vbs|\g(x,t;-v) \mbox{ converge  to a point in } S(f)\}.
$$
This is a compact set, and the  condition $(\gC)$
implies that this
set is a subset of
 $D(\indg {n-s} ; -v)$.
Therefore there is a \nei~ $U$
of $Q_s$ in
$\vbs$
\sut~
$\stv (U)$
is in
$D_\d(\indl s; v)\cap V_a$
and this last set is in
$\Int \vasm$
(again by $(\gC)$).
It follows that the map
$\stv$
gives rise to a well-defined continuous map
\bq
\vbs/\vbsm
\to
\vas/\vasm  \lbl{f:vdown}
\end{equation}
The family of these maps (one map for each
$s:0\leq s\leq n-1$)
is the substitute for
"cellular approximation"
of $\stv$.
Technically it is more convenient for us
to  quotient out a part of $\vbs$
a bit larger than $\vbsm$.
Recall that $\phi_1$
is an ordered Morse function;
let
$\a_0<\a_1<\dots <\a_n$
be the ordering sequence
for
$\phi_1$
\wrt~which $v$ is $\d$-separated.
Let
\bqq
V_b^{(s)}=
\phi_1^{-1}([\a_0, \a_{s+1}]), \quad
\vvbs=
D_\d(\indl s; u_0)\cup \vvbsm
\end{equation*}
Then
$\vbs\sbs \vvbs$
and it is easy to
see that the map (\ref{f:vdown})
is factored through the map
\bq
\vflesh
:
\factor
\to
\factora
\end{equation}
The pair
$(\vvbs, \vvbsm)$
is  homotopy equivalent to
$(\vbs, \vbsm)$,
but it has the advantage that
$\factor$
splits as the wedge of spaces
corresponding to the descending discs.
Namely, for $p\in S_s(f)$
let
$\Sigma(p)=
D(p,v)/D(p,v)\cap \vvbsm$
and
$\Sigma_\d(p)=
D_\d(p,v)/D_\d(p,v)\cap \vvbsm$.
Then
$\Sigma(p)$
is homeomorphic to
a $s$-dimensional sphere,
the inclusion
$\Sigma(p)\rInto
\Sigma_\d(p)$
is a homotopy equivalence, and
$\factor$
is the wedge of all
$\Sigma_\d(p)$
with $p\in S_s(f)$.
The \ho~induced by $\vflesh$
in homology is called
{\it
homological gradient descent } and denoted
$\HH_s(-v)$.
The key property of the homological gradient descent
allowing in particular to compute the boundary operators
in the Novikov complex
in terms of
$\HH_s(-v)$
is the following:

\it

Let $N$ be an oriented \sma~ of $V_b$,
\sut~$N\sbs \vvbs$
and $N\sm\Int \vvbsm$
 is compact.
Then the manifold $N'=\stv (N)$ is in $\vvas$ and
$N'\sm \Int \vvasm$
is compact and the fundamental class of $N'$ modulo
$\vvasm$
equals 
$\HH_s(-v)\([N]\)$.

\rm

\subsection{Circle valued Morse maps}
\label{su:circle}

Now we can proceed to the circle-valued Morse maps.
First of all we introduce  the basic terminology.
Let $f:M\to S^1$ be such a map.
Assume that $f$ is primitive, that is
$f_*:H_1(M)\to H_1(S^1)=\ZZZ$ is epimorphic.
To simplify the notation we shall assume that
$1\in S^1$ is a regular value for $f$.
Denote $f^{-1}(1)$ by $V$.
Let $\CC:\bar M\to M$ be
the infinite cyclic covering,
associated to $f$, and $F:\bar M\to\RRR$ be a lifting of $f$.
Set
$$V_\a=F^{-1}(\a), \quad
W=F^{-1}\([0,1]\), \quad
V^-=F^{-1}\(]-\infty,1]\).$$

Thus the  cobordism $W$ is  the result of
cutting $M$ along $V$.
The structure group
of $\CC$
is isomorphic to $\ZZZ$ and we choose the generator $t$
of this group so as to satisfy  $V_\a t=V_{\a-1}$.
Denote $Wt^s$ by $W_s$; then $\bar M$
is the union
$\cup_{s\in\ZZZ} W_s$, the neighbor copies
$W_{s+1}$ and $W_s$
 intersecting by
$V_{-s}$.
For any $k\in\ZZZ$ the restriction of $\CC$ to $V_k$
is a diffeomorphism
$V_k\to V$.
Endow $M$ with an arbitrary riemannian metric and lift it to
a $t$-invariant riemannian metric on $\bar M$.
Now $W$ is a riemannian cobordism, and
$t^{-1}:\dow=V_0\to\daw=V_1$
is an isometry.
We shall say that $v$ satisfies condition
$(\gC')$
if the $(F\vert W)$-gradient $v$ satisfies the condition $(\gC)$
from \su~\ref{su:cc},
and, moreover, the Morse functions $\phi_0, \phi_1$
and their gradients
$u_0, u_1$
can be chosen so as to
satisfy
$\phi_0(xt)=\phi_1(x), t_*(u_1)=u_0$.
The set of
Kupka-Smale
$f$-gradients $v$ satisfying $(\gC')$
will be denoted by
$\GG_0(f)$.
The set $\GG_0(f)$
is $C^0$-open and dense in $\GG(f)$
(this is a version of the theorem
\ref{t:cc}, see
\cite{pafest}, \S 8).
Let $v\in\GG_0(f)$.
The condition $(\gC')$ provides a Morse function
$\phi_1$ on $V_1$
together with its gradient $u_1$,
and  a Morse function
$\phi_0=\phi_1\circ t^{-1}$ on $V_0$
together with its gradient $u_0=t_*(u_1)$.
For every $k\in\ZZZ$ we obtain also an ordered Morse function
$\phi_k=\phi_0\circ t^{k}:V_k\to\RRR$
 and a $\phi_k$-gradient
$u_k=(t^{-k})_*(u_0)$.

The universal covering
$\PP:\wi M\to M$
 factors through
$\CC$, that is,  there is a covering
$p:\wi M\to \bar M$
with structure group
$H=\Ker\xi$
\sut~$\CC\circ p=\PP$.
For a subset $X\sbs \bar M$
we shall denote
$p^{-1}(X)$
also by
$\wi X$.
Let $v\in \GG_0(f)$.
For $k\in\ZZZ$ we have the $H$-invariant filtrations
$\wi V_k^{\{\leq s\}}$
of
$\wi V_k$
and the corresponding equivariant version
of the homological gradient descent. That is
for every $k\in\ZZZ, s\in \NNN$
we have a continuous
$H$-equivariant map

\bq
\vvfl:
\fac\to
\facm  \lbl{f:eqvdown}
\end{equation}
and the \ho~
$\wi \HH(-v)$, induced by
$\vvfl$
in homology.
Let us introduce the group
\bq
T_s=
\bigoplus_{k\in\ZZZ}
H_*(\fac) \lbl{f:ts}
\end{equation}

This abelian group
has the obvious structure of free $\ZZZ G$-module, and the
homomorphism
$(\vflesh)_*$
is a $\ZZZ G$-\ho~
of this module.
Let $\tau_s$
be its matrix.
Note that the matrix entries of
$\tau_s$
are in $\ZZZ G_{(-1)}$.
Thus the matrix
$1-\tau_s$
is invertible
over the ring $\Lxi$
and the image of
$1-\tau_s$
in $K_1(\Lxi)$ belongs to $\wh W$.
In the statement of the theorem below
we shall keep the notation
$1-\tau_s$
for the image
of this element in
$\wh W$.

\beth
[\cite{pafest}, Corollary 7.16]
\lb{t:tors}
Let $v\in \GG_0 (f)$.
There is
a chain homotopy equivalence
\begin{equation}
\phi: \wi C_*(v) \to
\cmxi \lbl{f:heq}
\end{equation}
\sut~
\begin{equation}
\tau(\phi )=
\prod_{s=0}^{n-1} (1-\tau_s)^{(-1)^{s}}\hfill
\lbl{f:tors}
\end{equation}
     \enth
$\qs$

In \cite{pafest}
this theorem is stated for the case of
the  maximal abelian covering
instead of universal covering, but
the  arguments are carried  over to the
present case without difficulty.

\bere
We take here the opportunity to note that the cited formula
from \cite{pafest}, and, consequently, the formula (3)
of the main theorem of \cite{pafest} contain a sign error. 
With the sign convention of \cite{milnWT}, accepted in
\cite{pafest},
the formula (3) of the main theorem of \cite{pafest}
should read:
$\tau(\phi|G)= -\ove{\zeta_L(-v)}$
instead of
$\tau(\phi|G)= \ove{\zeta_L(-v)}$
in the paper.
\enre

\section{Proof of the main theorem }
\label{s:prfclo}

Recall from \ref{su:circle}
the subset
$\GG_0(f)\sbs \GG(f)$.
Let $v\in \GG_0(f)$.
In order to prove the main theorem
it remains to show that the \ho
~$\gL$
evaluated on the right hand side of
(\ref{f:tors})
equals to
$-\eta_L(-v)$.
The \ma~$V$
being identified with $V_0\sbs \bar M$
inherits from $V_0$
the handle-like filtrations
$V^{[s]}, V^{(s)}$.
Let
$Cl^{[s]}(-v)$
be the subset of all the closed orbits of $(-v)$
passing through a point of
$V^{\{\leq s\}}\sm V^{\{\leq s-1\}}$.
It follows from the property $(\gC' )$
that the set $Cl(-v)$
is the disjoint union of its subsets
$Cl^{[s]}(-v)$.
Therefore if we define
 \begin{equation}
\eta_s(-v)=
\sum_{\g\in Cl^{[s]}(-v)} \ve(\g) \frac {\{ \g\}}{m(\g)}, \lbl{f:etas}
 \end{equation}
it suffices to prove that for every $s$ we have:

 \begin{equation}
(-1)^{s+1} \gL (1- \tau_s)
=
\eta_s(-v)\lbl{f:htau}
 \end{equation}

So we fix $s$ up to the end of the proof.
We shall abbreviate
$Cl^{[s]}(-v)$
to
$Cl$.
Recall
from \ref{su:circle}
that for every $k\in\ZZZ$
we have a continuous map
 \begin{equation}
\vvfl:
\fac
\to
\facm \end{equation}
constructed from the gradient descent map.

Let $B_k=\vvks/\vvksm$
and
$B=\sqcup_k B_k$.
We shall denote by $\o_k$
the point
$[\vvksm]$
 of $B_k$.
The point $\o_0$ will be also denoted by $*$.
There is a natural action of $G$ on $B$,
derived from the action
of
$G$ on $\wi M$;
every element
$g\in G_{(-n)}$
sends $B_k$
homeomorphically to
$B_{k-n}$
and
$\o_k$ to $\o_{k-n}$.
In particular every $B_k$
is $H$-invariant.
The map
$     \vvfl: B_k\to B_{k-1}$
is $H$-equivariant and
$\tau_s=   ( \vvfl)_*: \oplus_k H_s(B_k)  \to \oplus_k H_s(B_{k-1})$.

Now we shall make an additional auxiliary choice: we choose
and fix for every
critical point $p\in S_s(f)$
its lifting to $\wi V_0$.
The system of these liftings will be denoted by
$\Upsilon$.
The reason for introducing $\Upsilon$ into the game is the following:
the series $\eta_s(-v)$
is well defined only
in $\gG$
which is a quotient of
$P=\QQQ H_\rho [[t]]$.
The choice of
$\Upsilon$ enables us to
construct a lifting
of these power series to $P$.
On the other hand the choice of $\Upsilon$
determines a base in the free
$\ZZZ G$-module $T_s
=
\oplus_k H_s(B_k)  
$ from (\ref{f:ts})
and so we obtain also a representative for
$\gL(1-\tau_s)$
in $P$.
Finally we shall compare the two resulting elements of $P$ by a simple
application of Lefschetz-Dold fixed point formula.

Note first of all that
$\Upsilon$
determines liftings of the discs $D_\d(p,v)$ to $\wi V_0\sbs\wi M$,
and liftings of the thickened spheres
$\Sigma_\d(p)$
to $B_0$.
The union of these lifted spheres  is
 a subspace $\beta\sbs B_0$, homeomorphic to
the wedge of the thickened spheres themselves.
Now
\bqq
B_k=\vee_{g\in G_{(k)}}\beta\cdot g, 
\quad
B=\cup_{g\in G} \b\cdot g.
\end{equation*}
The \hog
~classes in $B_0$
of the liftings of the spheres  $\Sigma_p$
form
a $\ZZZ$-base in $H_s(\b)$
and a $\ZZZ G$-base in $H_s(B)$.
The  element of this base corresponding to a critical
point
$p\in S_s(f)$
will be denoted by
$[p]$.
Now  we can construct from $\Upsilon$ a  lifting of
$\gL(1-\tau_s)\in \bar P$
to $P=\QQQ H_\rho [[t]]$.
Let
$(R(k)_{qp})$
be the matrix of the
$\ZZZ G$-homomorphism
$\tau^k_s$
with respect to our  base.
That is
$\tau_s^k([p])=
\sum_{q\in S_s(f)} [q]\cdot R(k)_{qp}$
with
$R(k)_{qp}\in \ZZZ G_{(-k)}$.
Set
$\k_k=\sum_p R(k)_{pp}\in \ZZZ G_{(-k)}$
then the element
$$K(v)= -\sum_{k=1}^\infty \frac {\k_k}k$$
is obviously a lifting to
$P=\QQQ H_\rho  [[t]]$
of
$\gL(1-\tau_s)$.

To lift to $P$ the right hand side of (\ref{f:etas}) we shall first
 translate our data to the language of
 fixed point theory.
We shall say that a point $a\in\beta\sm\{*\}$
is a {\it $G$-fixed-point} of $(\vvfl)^k$, if
$(\vvfl)^k(a)=a\cdot g$.
The element $g\in G_{(-k)}$
is uniquely determined by $a$ and will be denoted by $g(a)$.
The set of all
$G$-fixed points  of $(\vvfl)^k$
will be denoted by $GF(k)$.
The set of all
$G$-fixed points  of $(\vvfl)^k$
with given $g(a)=g$
will be denoted by $GF(k,g)$.
Thus
$GF(k)=
\sqcup_{g\in G_{(-k)}}GF(k,g)$.
By analogy with standard fixed point theory we  define the
multiplicity
$\mu(a)$ and the index $\ind a$ for every $a\in GF(k)$.

Let $a\in GF(k)$. Let $a_i=(\vvfl)^i(a)$
and let $\wh a_i$ be (the
unique) point in $\b$ belonging to the $G$-orbit of $a_i$.
The set of all $\wh a_i, i\in\NNN$
will be called the {\it quasiorbit of $a$}, and denoted by
$\QQ (a)$; it is a finite subset of $\b$ of 
cardinality $\frac k{\mu(a)}$.
Consider the integral curve
$\g$ of $(-v)$ in $\wi M$, such that $\g(0)=a$; then for some
$T>0$ we have
$$
\g(T)=
(\vvfl)^k(a)=a\cdot g(a)\in \wi V_{-k}
$$
The map
$\PP\circ\g: [0,T]\to M$
 is then a closed orbit. Thus we obtain a map
$\a:       GF(k)\to Cl$
whose image is exactly the subset
$Cl_k\sbs Cl$, consisting of all
$\g\in Cl$
with
$f_*([\g])=-k$.
For every $\g\in Cl_k$
the set $\a^{-1}(\g)$
is the quasiorbit
of some $a\in GF_k$.
(Note that the set
$\PP(\a^{-1}(\g))$
is the intersection of the orbit $\g$ with $V$.)
Further, for every
$a\in \a^{-1}(\g)$
the projection of $g(a)$ to
the set $\G$
of conjugacy classes of $G$
equals to
$\{\g\}$.
Moreover,
$\ve(\g)=\ind a$
and
$m(\g)=\m(a)$.
 Thus the following power series
 \begin{equation}
 \nu(v) =
 \sum_{k=1}^\infty \frac 1k \sum_{a\in GF_k}
(\ind a) g(a)
 \label{f:nunu}
 \end{equation}
is a lifting to
$P= \QQQ H_\rho [[t]]$ of
$\eta_s(-v)$.

To prove our theorem it remains to show that
$\nu(v)=(-1)^{s+1} K(v)$. This follows obviously from the next lemma.

\bele
\lb{l:kappa}
\begin{equation}
\varkappa_k=(-1)^s\sum_{a\in GF(k)}  (\ind a) g(a)
\end{equation}
\enle
\Prf
The set $GF(k,g)$
is exactly the fixed point set of the following composition:
\begin{equation}
\beta=B_0\rTo^{\displaystyle{(\vvfl})^k} B_{(-k)}= 
\bigvee_{h\in G_{(-k)}}\beta\cdot h
\rTo^{\displaystyle{\pi_g}}\beta
\end{equation}
where $\pi_g$ is the map which sends every component of the wedge
to $*$, except the component $\beta g$, and this one is sent to
$\beta$ via the map $x\mapsto x g^{-1}$.
The point $*\in\b$
is a fixed point of this map, and its index equals to 1.
Applying the Lefschetz-Dold fixed point formula,
the proofs of the Lemma and the main theorem are now complete.
$\qs$

\section{Further remarks}
\label{s:further}

In this section we shall discuss some particular cases of our
results and open questions.

\subsection{Fibrations over a circle}
\label{su:fibr}
Consider the case of maps $f:M\to S^1$
without critical points. The main theorem
then implies that the torsion of
the acyclic
$\Lxi$-complex
$
\cmxi=
C_*^\D(\wi M)\tens{\L}\Lxi $
satisfies
\bq\label{f:no_crit}
\gL(\tau(\cmxi))=-\eta_L(-v)
\end{equation}

The gradient descent in this case gives rise to an everywhere
defined diffeomorphism, say,
$g$
of $V=f^{-1}(\l)$
to itself (where $\l$ is a regular value
of $f$).
The closed orbits of $v$ are in one-to-one
correspondence with
periodic points
of $g$, and the formula
(\ref{f:no_crit})
thus gives an expression of a certain power series
obtained
via counting periodic points
(eta function of $g$)
in homotopy invariant terms.
It would be interesting to compare this formula
with other non-abelian invariants
counting periodic points of diffeomorphisms
like the invariants of Jiang (see \cite{jiang})
and Geoghegan-Nicas
(see \cite{geo}).

\subsection{Case of abelian fundamental group}
\label{su:abel}
Let us assume now that $G=\pi_1(M)$
is abelian, so that
$G=H\times \ZZZ$
and
$\xi:G\to\ZZZ$
is the projection onto the second direct summand.
In this case
the map
$W\to \wh W$
is an isomorphism
(since for any  commutative ring $R$ we have
$K_1(R)\approx R^\bu\oplus SK_1(R)$).
The element
$\tau(\phi)$
from the main theorem
can therefore be identified
with an element
of $W$.
Another attractive feature of the commutative case
is that
the
exponent of the element $\eta(-v)$
is well defined
and belongs to the Novikov ring.
This exponent is called {\it homological zeta function}
and denoted by $\zeta_L(-v)$.
The main theorem can now be reformulated
as follows:
\bq\label{f:ab_case}
\tau(\phi) =
(\zeta_L(-v))^{-1}
\end{equation}
where both sides are elements of the group $W$.
This formula is proved in \cite{pafest}.
Actually one can say more: both the sides of
(\ref{f:ab_case}) are rational functions (and not merely power series).
See \cite{pafest} for details.

\subsection{The simplest possible example}
\label{su:example}

Let $M=S^1$, and
consider the identity map
$f=\id:S^1\to S^1$.
Then $\bar M=\wi M=\RRR$,
choose the identity map as the lift of $f$ to $\RRR$.
Let $t$ be the generator of $\ZZZ$
acting on $\RRR$ as follows:
$x\mapsto x-1$.
For every $k\geq 1$
there is a unique orbit
of $(-v)$
in the class $t^k$;
this orbit
has index $1$.
Thus
$\eta_L(-v)=
\sum_{k\geq 1} t^k/k$.
The
$\ZZZ[t,t^{-1}]$-complex
$C_*(\RRR)$
is of the form
\bq
\{0\lTo \L\lTo^\pr \L \lTo 0 \}
\end{equation}
where $\pr(1)=1-t$.
Tensoring over $\wh\L=\ZZZ((t))$
makes this complex
acyclic with torsion
$\tau=\tau(C_*(\RRR)\tens{\L}\wh \L)=1-t$
(see \cite{milnWT}, p.387)
Thus
\bqq
\ln\tau =\ln(1-t)=-\sum_{k\geq 1}   \displaystyle{\frac {t^k}k} 
=-\eta(-v)
\end{equation*}

\subsection{About Witt vectors}
\label{su:about_witt}

The theorem \ref{t:kerl}
together with Corollary \ref{c:laur_ser}
implies  that $\ove{\Log}:W/\ove{W_0} \to \bar P_+$
is an isomorphism,
where 
$W_0= \Ker(W\to K_1(R))$.
We do not know whether
$\ove{W_0}\sbs
\Ker(W\to K_1(P)$).
\pa
\unde{ Problem}\quad
Compute
$\ove{W_0}/W_0$.
This abelian group is an invariant of the ring $A$
and its automorphism $\rho$,
which will be denoted by
 $\gW(A, \rho)$.
Is it true that
 $\gW(A,\rho)$
vanishes in simple cases? for example when 
$A=\QQQ H$ where $H=\Ker (\xi:G\to\ZZZ)$ is a finite group?
(when $G$ is abelian  the
group
 $\gW(\QQQ H, \id)$
vanishes by the obvious reasons).

This context suggests that the natural power series topology in the
Novikov rings may have some deeper geometrical  sense.
It would be interesting to investigate
the "continuous Whitehead groups"
of the Novikov rings, defined by
\bqq
K_1^c(\Lxi)=
\GL(\infty, \Lxi) \Big/ \ove{  \EE  }
\end{equation*}
(where
$\EE=
[\GL(\infty, \Lxi), \GL(\infty, \Lxi)]  $
is the commutator of the general linear group, and 
$
\ove{ \EE}
$
is its closure ).

\subsection{Two questions about gradient flows}
\label{su:two_quest}
\pa
{\bf 1.} \quad 
Note that both parts of the equality
(\ref{f:ntzeta2})
are defined for every 
Kupka-Smale $f$-gradient $v$, not 
only for $C^0$-generic one. 
\pa
\unde{Question:} \quad 
Does the equality 
(\ref{f:ntzeta2})
hold also for arbitrary Kupka-Smale gradients?
\pa
{\bf 2.}
Assume now that $\eta_L(-v)=0$.
Then the element $\tau(\phi)$,
which is in the subgroup
$\wi W\sbs K_1(R)$
by (\ref{f:tors}),
belongs to 
$\Ker\gL:K_1(R) \to \bar P_+$.
Therefore, the element
$\tau(\phi)$
belongs to the subgroup
$\ove{W_0}/W_0\sbs K_1(R)$
where $W_0=\Ker J$ and $\ove{W_0}$ is the closure of the subgroup $W_0$
in the natural topology, see \ref{s:applic_log} and
\ref{su:about_witt}.
\pa
\unde{Question:} \quad 
What is the geometric meaning of the image of $\tau(\phi)$
in 
$\gW(\QQQ H, \rho)= \ove{W_0}/W_0$?
Can this element be interpreted in terms of closed orbits of
the flow?

\subsection{Constructions with non-commutative localization}
\label{su:const_local}

One can show
(see \cite{paepri})
that for $C^0$-generic gradients
the Novikov
complex
$C_*(v)$
can be defined
over a non-commutative localization
$\L_{(\xi)}$
of the group ring $\L=\ZZZ\pi_1(M)$.
In the paper \cite{farran}
M.Farber and A.Ranicki
give a purely algebraic construction
of a chain complex
$C_*^{FR}(M,f)$
over $\L_{(\xi)}$
built from Morse-theoretic data associated to
a circle-valued Morse function $f:M\to  S^1$
and
counting the localized homology of the universal covering.
In the paper
\cite{rann}
A.Ranicki constructs an explicit
isomorphism
between the Novikov complex and the complex $ C_*^{FR}$.
It would be interesting to compare this isomorphism
with the chain equivalence $\phi$ constructed in the present paper.

\label{refer}

\end{document}